\documentclass[11pt, reqno]{amsart}

\usepackage{amssymb,amsthm,amsmath,amsfonts}
\usepackage[margin=1in]{geometry}

\usepackage{mathrsfs}
\usepackage{verbatim}
\usepackage{color}
\usepackage{enumitem}
\usepackage{bm}
\usepackage[parfill]{parskip}
\usepackage[usenames,dvipsnames]{xcolor}
\usepackage[colorlinks=true, linkcolor=blue, citecolor=red]{hyperref}
\usepackage[hyphenbreaks]{breakurl}

\def\mathcal{\mathscr}

\theoremstyle{plain}
\newtheorem{theorem}{Theorem}[section]

\theoremstyle{remark}
\newtheorem{remark}[theorem]{Remark}

\theoremstyle{plain}
\newtheorem{corollary}[theorem]{Corollary}

\newtheorem{proposition}[theorem]{Proposition}
\newtheorem{definition}[theorem]{Definition}

\numberwithin{equation}{section}

\makeatletter
\newcommand{\LeftEqNo}{\let\veqno\@@leqno}
\makeatother

\renewcommand{\d}{\mathrm{d}}
\DeclareMathOperator*{\R}{{\mathbb{R}}}
\DeclareMathOperator*{\E}{\,{\mathbb{E}}}
\renewcommand*{\P}{\,\mathbb{P}}
\DeclareMathOperator*{\Pas}{\,\mathbb{P}\,-\,\mathrm{a.s.}}

\DeclareMathOperator*{\Dom}{\mathrm{Dom}}
\renewcommand{\L}{\mathcal{L}}
\renewcommand{\O}{\mathcal{O}}

\begin{document}
\title[]
{$\bm{L^p}$-valued stochastic convolution integral driven by Volterra noise}
\author{P. \v{C}oupek}
\address{Charles University\\
Faculty of Mathematics and Physics\\
Sokolovsk\'{a} 83\\
186 75\\
Prague 8\\
Czech Republic}
\email{coupek@karlin.mff.cuni.cz}

\author{B. Maslowski}
\address{Charles University\\
Faculty of Mathematics and Physics\\
Sokolovsk\'{a} 83\\
186 75\\
Prague 8\\
Czech Republic}
\email{maslow@karlin.mff.cuni.cz}

\author{M. Ondrej\'{a}t}
\address{Czech Academy of Sciences\\
Institute of Information Theory and Automation\\
Pod Vod\'{a}renskou v\v{e}\v{z}\'{i} 4\\
182 08\\
Prague 8\\
Czech Republic}
\email{ondrejat@utia.cas.cz}

\thanks{The first author was supported by the Charles University, project GAUK No.\ 322715 and SVV 2016 No.\ 260334. The second and the third authors were supported by the Czech Science Foundation, project GA\v{C}R No.\ 15-08819S.\\
\indent Electronic version of an author accepted manuscript of an article published as P.~\v{C}oupek, B.~Maslowski, and M.~Ondrej\'at, \textit{$L^p$-valued stochastic convolution integral driven by Volterra noise}, Stochastics and Dynamics, Vol.\ 18, No.\ 6, 1850048 (2018) DOI: \href{https://www.worldscientific.com/doi/abs/10.1142/S021949371850048X}{10.1142/S021949371850048X} {\textcopyright} World Scientific Publishing Company, \url{http://www.worldscientific.com/worldscinet/sd}.}

\begin{abstract}
Space-time regularity of linear stochastic partial differential equations is studied. The solution is defined in the mild sense in the state space $L^p$. The corresponding regularity is obtained by showing that the stochastic convolution integrals are H\"{o}lder continuous in a suitable function space. In particular cases, this allows to show space-time H\"{o}lder continuity of the solution. The main tool used is a hypercontractivity result on Banach-space valued random variables in a finite Wiener chaos.
\end{abstract}
\maketitle

\section{Introduction}

The paper is devoted to the study of mild solutions to linear stochastic differential equations in the Lebesgue $L^p(D,\mu)$ space perturbed by additive noise of Volterra type. Sufficient conditions for the existence and H\"{o}lder continuity of the solutions are given. This allows to show that solutions to particular stochastic partial differential equations (SPDEs) are H\"{o}lder continuous random fields. 

More precisely, we consider the stochastic evolution equation which takes the form
	\begin{equation}
	\label{eq:SDE_0}
		\d{X}_t = AX_t+\varPhi\,\d{B}_t, \quad X_0=x,
	\end{equation}
where $B$ is an infinite-dimensional $\alpha$-regular Volterra process which belongs to a finite Wiener chaos and $A$ is a generator of a strongly continuous, analytic semigroup $(S(t), t\geq 0)$ of operators acting on the space $L^p(D,\mu)$ with $1<\frac{2}{1+2\alpha}\leq p<\infty$. The mild solution is given by the stochastic convolution integral
	\begin{equation*}
		X_t = S(t)x +\int_0^tS(t-r)\varPhi\d{B}_r, \quad t\geq 0,
	\end{equation*}
and we give sufficient conditions for its existence and H\"{o}lder continuity in the domain of a fractional power of $A$. Canonical examples of SPDEs to which our theory may be applied are the heat equation on bounded domain $\O\subset\R^d$ with pointwise noise, formally given by 
	\begin{equation*}
		\partial_t u(t,\xi) = \Delta u (t,\xi) + \eta(t)\delta_z(\xi) \quad\mbox { on }\quad {\R}_{+}\times \mathcal{O},
	\end{equation*}
where $\delta_z$ is the Dirac distribution; or the heat equation with distributed noise
	\begin{equation*}
		\partial_t u(t,\xi) = \Delta u(t,\xi) + \eta(t,\xi) \quad\mbox { on }\quad  {\R}_+\times\O.
	\end{equation*}
In both these examples, our results allow to make use of the embdedding of Sobolev-Slobodeckii spaces into the spaces of H\"{o}lder continuous functions and hence, by taking sufficiently large $p$, to obtain space-time H\"{o}lder continuity for $d=1,2,3$. 

The scalar Volterra process is a stochastic process which might not be Markov, Gaussian or a semimartingale but which admits a certain covariance structure instead. In particular, there is a kernel $K$ which satisfies suitable regularity conditions (see subsection \ref{sec:VP}) such that the covariance function can be written as 
	\begin{equation*}
		R(s,t) = \int_0^{s\wedge t}K(s,r)K(t,r)\d{r}, \quad s,t\geq 0.
	\end{equation*}
The most notable examples which satisfy the definition are the fractional Brownian motion (fBm) with the Hurst parameter $H>1/2$ which is Gaussian and lives in the first Wiener chaos (see e.g. \cite{AlosNua03, DecrUstu99, MandelbrotVanNess} for its definition and further properties), and the Rosenblatt process, which is non-Gaussian and lives in the second Wiener chaos (see e.g. \cite{Taqqu11, Tud08}). Note, however, that the class of $\alpha$-regular Volterra processes in a finite Wiener chaos is not restricted only to these two processes (see \cite{CouMas16} for other examples). See also \cite{AlosMazNua01, BauNua03,BrzNeeSal12, ErrEss09, Hida60, Hult03}.

Stochastic convolution integral with respect to $\alpha$-regular Volterra processes has already been considered in \cite{CouMas16} in the Hilbert space setting. In particular, it has been shown that the integral admits a version with H\"{o}lder continuous sample paths of very small order which can be improved by $1/2$ if the driving process is Gaussian. In the present paper, we further develop this idea by assuming that the Volterra process lives in a finite Wiener chaos which allows us to prove hypercontractivity of the $L^p$-valued stochastic integrals. This yields the same regularity in the non-Gaussian case as in the Gaussian case. Other works on evolution equations driven by Volterra processes are \cite{CouMasSnup17, DunMasDun17}. See also \cite{TudBon11} where the authors consider stochastic convolution in Hilbert spaces driven by processes from a finite Wiener chaos which have similar covariance structure as $\alpha$-regular Volterra processes.

The paper is organized as follows.

In Section 2, we collect the tools needed in the following sections. In particular, definition of an $\alpha$-regular Volterra process is given and one-dimensional stochastic integration of deterministic real-valued functions is recalled. This part closely follows the papers \cite{AlosMazNua01, CouMas16}. We then proceed to the definition of the $n$-th Wiener chaos and give a hypercontractivity result which states that Banach-space valued linear combinations of elements in $n$-th Wiener chaos have equivalent moments (Proposition \ref{prop:HC_in_chaos}). Then, definition of $\gamma$-radonifying operators follows together with their basic properties. Finally, we collect the main assumptions used in the paper. 

Section 3 is devoted to stochastic integration of operator-valued functions with respect to cylindrical Volterra process in the space $L^p(D,\mu)$. We give characterization of admissible integrands (Proposition \ref{prop:charakterization_of_integrands}) and identify sufficient conditions for integrability on the scale of Lebesgue spaces (Corollary \ref{cor:bound_elem_Lp_int}). 

Section 4 contains the main results of the paper. In particular, sample path measurability of the mild solution to \eqref{eq:SDE_0} is shown under certain natural conditions on the semigroup $S$ (Proposition \ref{prop:meas_paths}) and then, factorization method is used to prove H\"{o}lder continuity of the solution under slightly stronger conditions (Proposition \ref{prop:V_dp-continuity}).

The paper is concluded with two examples contained in Section 5 - the stochastic heat equation with pointwise Volterra noise and the stochastic parabolic equation of the $2m$-th order with distributed noise which is Volterra in time and can be both white or correlated in space.

\section{Preliminaries}

This section collects the main tools needed in the next sections. Throughout the paper, $A\lesssim B$ means that there exists a positive constant $C$ such that $A\leq CB$. Similarly, the symbol $A\eqsim B$ means that there exist positive constants $C_1, C_2$ such that $C_1B\leq A\leq C_2 B$.

\subsection{Volterra processes}
\label{sec:VP}

Consider a measurable function $K: {\R}_{+}\times {\R}_{+}\rightarrow {\R}_{+}$ which is
\begin{itemize}
	\item Volterra, i.e.
		\begin{enumerate}
				\item[(i)] $K(0,0)=0$ and $K(t,r)=0$ on $\{0\leq t < r<\infty\}$,
				\item[(ii)] $\lim_{t\rightarrow r+}K(t,r) = 0$ for all $r\geq 0$;
		\end{enumerate}
	\item and $\alpha$-regular, i.e.
		\begin{enumerate}
				\item[(iii)] $K(\cdot,r) \in\mathcal{C}^1(r,T)$ for all $r\geq 0$, $T>0$; and there is an $\alpha\in (0,1/2)$ such that
						\begin{equation*}
							\left|\frac{\partial K}{\partial u}(u,r)\right| \lesssim (u-r)^{\alpha-1}\left(\frac{u}{r}\right)^\alpha
						\end{equation*}
						on $\{0<r<u\leq T\}$.
		\end{enumerate}
\end{itemize}
Such a function $K$ is called an \textit{$\alpha$-regular Volterra kernel} in the sequel.

\begin{definition}
\label{def:VP}
We say that a real, centered stochastic process $b=(b_t, t\geq 0)$ defined on a probability space $(\Omega,\mathcal{F},\P)$ is an $\alpha$-regular Volterra process if $b_0=0$ $\Pas$, and its covariance function takes the form
	\begin{equation}
	\label{eq:covariance_function_of_VP}
		{\E}(b_sb_t) = R(s,t) = \int_0^{s\wedge t}K(s,r)K(t,r)\mathrm{d}r,  \quad s,t\geq 0,
	\end{equation}
for some $\alpha$-regular Volterra kernel $K$.
\end{definition}

\begin{remark}
Note that if $K$ is an $\alpha$-regular Volterra kernel, it holds that 
	\begin{equation}
	\label{eq:continuity}
		\int_0^t(K(t,r)-K(s,r))^2\d{r} \lesssim (t-s)^{1+2\alpha}, \quad 0\leq s<t,
	\end{equation}
using the fact that
	\begin{equation*}
		(uv)^\alpha\int_0^{u\wedge v} r^{-2\alpha}(u-r)^{\alpha-1}(v-r)^{\alpha-1}\d{r} = \mathrm{B}(\alpha,1-2\alpha)|u-v|^{2\alpha-1}
	\end{equation*}
holds for $u,v\geq 0, u\neq v$, where $\mathrm{B}$ is the Beta function. It follows, in particular, that $K(t, \cdot)\in L^2(0,t)$ for every $t>0$ which makes the integral on the right-hand side of \eqref{eq:covariance_function_of_VP} finite for every $s,t\geq 0$. Moreover, using \eqref{eq:covariance_function_of_VP}, \eqref{eq:continuity} and the Kolmogorov continuity criterion, we can infer that $b$ has a version with $\varepsilon$-H\"older continuous sample paths for every $\varepsilon\in (0,\alpha)$. If, additionally, the process $b$ is assumed to live in a finite Wiener chaos (see the forthcoming subsections \ref{sec:chaos} and \ref{sec:hypotheses}), its increments have equivalent moments (see Proposition \ref{prop:HC_in_chaos}) and the Kolmogorov continuity criterion implies that $b$ has a version with $\varepsilon$-H\"older trajectories for every $\varepsilon\in (0,\alpha+1/2)$.
\end{remark}

\begin{remark}
The fractional Brownian motion (fBm) with the Hurst parameter $H\in (1/2,1)$ is an $\alpha$-regular Volterra process with $\alpha=H-1/2$. The fBm is defined as the centered Gaussian process with continuous sample paths whose covariance function is 
	\begin{equation*}
		R^H(s,t) := \frac{1}{2}\left(|s|^{2H} + |t|^{2H} - |t-s|^{2H}\right), \quad s,t\geq 0.
	\end{equation*}
The function $R^H(s,t)$ can indeed be written as \eqref{eq:covariance_function_of_VP} with
	\begin{equation}	
	\label{eq:kernel_of_fBm}
		K^H(t,r) = C_H\int_r^t\left(\frac{u}{r}\right)^{H-\frac{1}{2}}(u-r)^{H-\frac{3}{2}}\,\mathrm{d}u, 
	\end{equation}
where $C_H$ is a suitable constant (see \cite{AlosNua03}). Although the fBm is Gaussian, in the present paper, Gaussianity is not assumed. The Rosenblatt process $Z$ is an example of a non-Gaussian Volterra process. In particular, $Z$ is defined as
	\begin{equation*}
		Z_t := C\int_{\R}\int_{\R}\left(\int_0^t(u-y_1)_{+}^{-\frac{2-H'}{2}}(u-y_2)_{+}^{-\frac{2-H'}{2}}\d{u}\right)\d{W}_{y_1}\d{W}_{y_2}, \quad t\geq 0,
	\end{equation*}
where $C$ is a normalizing constant such that ${\E}Z_1^2=1$, $H'\in (1/2,1)$ and $W$ is the two-sided standard Wiener process. Then, $Z$ satisfies Definition \ref{def:VP} with the kernel $K^{H'}$ which takes the form \eqref{eq:kernel_of_fBm} and hence, $Z$ is an $\alpha$-regular Volterra process with $\alpha=H'-1/2$ (see \cite{Taqqu11} and \cite{Tud08}). Other examples of $\alpha$-regular Volterra processes include the Liouville fractional Brownian motion (see \cite{BrzNeeSal12}) or the Liouville multifractional Brownian motion (LmfBm) where the Hurst parameter $H=H(t)$ may be a function of $t$. Its definition and the assumptions on $H(\cdot)$ which ensure that the LmfBm is an $\alpha$-regular Volterra process are given in \cite[Example~2.14]{CouMas16}.
\end{remark}

\subsection{Wiener integral}
Since Volterra processes are not necessarily semimartingales, the standard It\^{o} approach to a stochastic integral is not applicable. The (already rather standard) definition of Wiener-type integrals (i.e. for deterministic integrands) driven by scalar Volterra processes is given below (cf. \cite{AlosMazNua01} and \cite{CouMas16}).

Let $T>0$ and consider the linear space of (${\R}$-valued) deterministic step functions $\mathcal{E}$, i.e.
	\begin{align*}
		\mathcal{E} & := \bigg\{\varphi:[0,T]\rightarrow \R, \varphi = \sum_{i=1}^{n-1}\varphi_i\bm{1}_{[t_{i},t_{i+1})}+\varphi_{n}\bm{1}_{[t_n,t_{n+1}]},\\
		& \hspace{2cm}\varphi_i\in \R, i\in\{1,\cdots,n\}, 0=t_1< t_2<,\cdots,< t_{n+1}=T, n\in\mathbb{N}\bigg\}.
	\end{align*}
Define an operator $\mathcal{K}_T^*:\mathcal{E}\rightarrow L^2(0,T)$ by
	\begin{equation}
	\label{eq:Kstar}
			\big(\mathcal{K}^*_T\varphi\big) (r) := \int_r^T\varphi(u)\frac{\partial K}{\partial u}(u,r)\d{u}, \quad \varphi\in\mathcal{E}.
		\end{equation}
Let $b=(b_t, t\geq 0)$ be an $\alpha$-regular Volterra process with the kernel $K$. Consider the linear operator $i_T:\mathcal{E}\rightarrow L^2(\Omega)$ given by
	\begin{equation*}
		\varphi := \sum_{i=1}^{n-1}\varphi_i\bm{1}_{[t_{i},t_{i+1})}+\varphi_{n}\bm{1}_{[t_n,t_{n+1}]}\quad\overset{i_T}{\longmapsto}\quad \sum_{i=1}^{n-1}\varphi_i(b_{t_{i+1}}-b_{t_i}) =:i_T(\varphi).
	\end{equation*}
Using \eqref{eq:covariance_function_of_VP} and \eqref{eq:Kstar}, it can be shown that
	\begin{equation}
	\label{eq:Itoisometry}
		\|i_T(\varphi)\|_{L^2(\Omega)} = \|\mathcal{K}^*_T\varphi\|_{L^2(0,T)}
	\end{equation}
which is an It\^{o}-type isometry for Volterra processes. For $f,g\in\mathcal{E}$, set 
\begin{equation*}
	\langle f,g\rangle_\mathcal{D}:=\langle\mathcal{K}^*_Tf,\mathcal{K}^*_Tg\rangle_{L^2(0,T)}.
\end{equation*} Without loss of generality, we assume that $\mathcal{K}^*_T$ is injective and thus, the function $\langle\cdot,\cdot\rangle_\mathcal{D}$ defines an inner product on $\mathcal{E}$. If this is not the case, we consider the quotient space $\tilde{\mathcal{E}}:=\mathcal{E}/\ker\mathcal{K}_T^*$ and we lift $\mathcal{K}^*_T$ to $\tilde{\mathcal{E}}$. Completing $\mathcal{E}$ under $\langle\cdot,\cdot\rangle_{\mathcal{D}}$	
yields the Hilbert space $(\mathcal{D},\langle\cdot,\cdot\rangle_{\mathcal{D}})$ and extends $\mathcal{K}^*_T$ to $\mathcal{D}$. This, in turn, extends $i_T$ to an operator from $\mathcal{D}$ to $L^2(\Omega)$ by \eqref{eq:Itoisometry}. $\mathcal{D}$ is the space of admissible integrands with respect to $b$ and $i_T:\mathcal{D}\rightarrow L^2(\Omega)$ is the Wiener-type integral. The usual notation is $i_T(\varphi) = \int_0^T\varphi\d{b}$. The space $\mathcal{D}$ can be very large and thus, it is important to identify certain function spaces which can be embedded into $\mathcal{D}$. By \cite[Proposition 2.9]{CouMas16} and its proof, we have that
	\begin{equation}
	\label{eq:estimate_norm_of_elementary_integral}
		\|i_T(\varphi)\|_{L^2(\Omega)}^2 \lesssim \int_0^T\int_0^T\varphi(u)\varphi(v)|u-v|^{2\alpha-1}\d{u}\d{v}
	\end{equation}
from which it follows that
	\begin{equation*}
		 L^\frac{2}{1+2\alpha}(0,T)\hookrightarrow \mathcal{D}.
	\end{equation*}

\subsection{Cylindrical Volterra processes}
In order to consider stochastic evolution equations, Vol\-ter\-ra processes with values in Hilbert spaces need to be introduced. The following definition provides such a generalization.

\begin{definition}
\label{def:cVP}
Let $(\Omega,\mathcal{F},\mathbb{P})$ be a probability space and let $U$ be a real separable Hilbert space. The \textit{$\alpha$-regular $U$-cylindrical Volterra process} is a collection $B=(B_t, t\geq 0)$ of bounded linear operators $B_t:U\rightarrow L^2(\Omega)$ such that
	\begin{itemize}
	\itemsep-0.5em
		\item for every $u\in U$, $B(u)$ is a real, centered stochastic process and $B_0(u)=0$ $\Pas$;
		\item for every $s,t\geq 0$ and every $u,v\in U$, it holds that
				\begin{equation*}
					\E B_s(u)B_t(v) = R(s,t)\langle u,v\rangle_U,
				\end{equation*}
				where $R(s,t)$ is given by \eqref{eq:covariance_function_of_VP} with some $\alpha$-regular Volterra kernel $K$.
	\end{itemize}
\end{definition}

\begin{remark}
\label{rem:representation} 
Let $\{e_n\}$ be a complete orthonormal basis of $U$. One may think of the process $B$ as the (formal) sum
	\begin{equation}
	\label{eq:representation}
		\tilde{B}_t = \sum_ne_nb_t^{(n)}
	\end{equation}
where $b_t^{(n)}:=B_t(e_n)$. The sequence $\{b^{(n)}\}$ consists of mutually uncorrelated (not necessarily independent) one-dimensional $\alpha$-regular Volterra processes. Similarly, as in the case of the standard cylindrical Brownian motion, the sum \eqref{eq:representation} does not converge in $L^2(\Omega;U)$. However, the integral with respect to $B$ introduced in section \ref{sec:stoch_integration_Lp} is well-defined as a random variable with values in a certain $L^p$-space.
\end{remark}

\subsection{Wiener chaos}
\label{sec:chaos}
Let further $\mathfrak{H}$ be a real separable Hilbert space and let $W$ be an $\mathfrak{H}$-isonormal Gaussian process, i.e. $W=(W(h), h\in \mathfrak{H})$ is a centered Gaussian family such that 
	\begin{equation*}
		\E W(h_1)W(h_2) := \langle h_1, h_2\rangle_\mathfrak{H}, \quad h_1, h_2\in \mathfrak{H}.
	\end{equation*}
Denote by $H_n$ the $n$-th Hermite polynomial, i.e.
	\begin{equation*}
		H_n(x) := \frac{(-1)^n}{n!}e^\frac{x^2}{2}\frac{\d^n}{\d{x}^n}\left(e^{-\frac{x^2}{2}}\right). 
	\end{equation*}
The $n$-th Wiener chaos (of $W$), $\mathcal{H}_n$, is the closed linear subspace of $L^2(\Omega)$ generated by the linear span $\{H_n(W(h)), h\in \mathfrak{H}, \|h\|_\mathfrak{H}=1\}$. In particular, the space $\mathcal{H}_0$ consists of constant random variables and $\mathcal{H}_1$ contains zero-mean Gaussian random variables which can be interpreted as stochastic integrals (with respect to $W$). For further reference see e.g. \cite{Nualart06}. We shall use the following feature of the spaces $\mathcal{H}_n$:

\begin{proposition}\label{prop:HC_in_chaos}
Let $p,q\in[1,\infty)$ and $n\ge 0$. Then there exists a number $C_{p,q,n}$ such that 
\begin{equation}\label{eq:chaos_moments}
\left\|\sum_{j=1}^m\xi_j x_j\right\|_{L^q(\Omega;X)}\leq C_{p,q,n}\left\|\sum_{j=1}^m\xi_j x_j\right\|_{L^p(\Omega;X)}
\end{equation}
holds for every Banach space $X$, $m\in\mathbb{N}$, and every $\{\xi_j\}_{j\le m}\subseteq\bigoplus_{i=0}^n\mathcal{H}_i$ and $\{x_j\}_{j\le m}\subseteq X$.
\end{proposition}

\begin{proof}
One can prove the inequality \eqref{eq:chaos_moments} either by a decoupling argument and the Kahane-Khintchine inequality as in \cite[Proposition 3.1]{maas} or by the scalar Neveu inequality \cite{neveu} applied to vector-valued Mehler's formula for the Ornstein-Uhlenbeck semigroup, see \cite[Theorem 1.4.1]{Nualart06} and the remark on page 62 in \cite{Nualart06}. The latter approach yields the inequality \eqref{eq:chaos_moments} only for $p,q\in(1,\infty)$ which then extends to the general range by the H\"older inequality, see the remark following \cite[Theorem 3.2.2, p. 113-114]{pena_gine}.
\end{proof}
	
\subsection{\texorpdfstring{$\bm{\gamma}$}{γ}-Radonifying operators} 
Let $U$ be a real separable Hilbert space and $E$ be a real separable Banach space. A bounded linear operator $R\in\L(U,E)$ is \textit{$\gamma$-radonifying} provided that there exists a centered Gaussian probability measure $\nu_R$ on $E$ such that
\begin{equation*}
	\int_E\varphi^2(x)\,\nu_R(\d{x})=\|R^*\varphi\|^2_U,\qquad\varphi\in E^*.
\end{equation*}
Such a measure is at most one, therefore we set
\begin{equation*}
	\|R\|^2_{\gamma(U,E)} :=\int_E\|x\|_E^2\,\nu_R(\d{x})
\end{equation*}
and denote by $\gamma(U,E)$ the space of $\gamma$-radonifying operators. It is well-known that $\gamma(U,E)$ equipped with the norm $\|\cdot\|_{\gamma(U,E)}$ is a separable Banach space (see \cite{Neidhardt78} or \cite{Ondrejat04}).

\begin{proposition}
\label{prop:radonifying_sum}
Let $R\in\L(U,E)$, $\{e_n\}$ be an orthonormal basis of $U$ and $\{\delta_n\}$ be a~sequence of independent standard centered Gaussian random variables. Denote
\begin{equation*}
	S_n:=\sum_{k=1}^n\delta_kRe_k.
\end{equation*}
The following claims are equivalent:
\begin{itemize}
\itemsep-0.5em
	\item $R$ is $\gamma$-radonifying,
	\item the sequence $\{S_n\}$ is convergent almost surely,
	\item the sequence $\{S_n\}$ is convergent in probability,
	\item the sequence $\{S_n\}$ is convergent in every $L^q(\Omega;E)$, $1\le q<\infty$.
\end{itemize}
\end{proposition}

\begin{proof}
This is a consequence of Theorem 2.3, Chapter V.2.4 and Theorem 5.3, Chapter V.5.3 in \cite{VakTarChob87}. Alternatively, the proof may be inferred from the It\^{o}-Nisio theorem and the Fernique theorem.
\end{proof}

\begin{remark}
\label{rem:separable}
Separability of a measure space $(D,\mu)$ means that there exists a countable system $\{V_n\}$ of measurable sets satisfying $\mu(V_n)<\infty$ such that, for every $\varepsilon>0$ and every measurable set $C$ satisfying $\mu(C)<\infty$, there is $n\in\mathbb{N}$ such that $\mu[(C\setminus V_n)\cup(V_n\setminus C)]<\varepsilon$. The following conditions are equivalent:
\begin{itemize}
\itemsep-0.5em
	\item the measure space $(D,\mu)$ is separable,
	\item there exists $1\le p<\infty$ such that $L^p(D,\mu)$ is a separable Banach space,
	\item $L^p(D,\mu)$ is a separable Banach space for every $1\le p<\infty$.
\end{itemize}
\end{remark}	

\begin{proposition}
\label{prop:radon_norm}
Let $(D,\mu)$ be a separable $\sigma$-finite measure space. Let further $1\leq p<\infty$ and $R\in\L(U, L^p(D,\mu))$. Then, $R\in\gamma(U, L^p(D,\mu))$ if and only if there exists a measurable function $r:D\to U$ satisfying
\begin{equation*}
	\int_D\|r(x)\|_U^p\,\mu(\d{x})<\infty,
\end{equation*}
and such that $(Ru)(x)=\langle r(x),u\rangle_U$ holds $\mu$-almost everywhere for every $u\in U$. There also exists a constant $c>1$, independent of $R$, such that
\begin{equation*}
	\frac{1}{c}\|R\|_{\gamma(U, L^p(D,\mu))}\le\left[\int_D\|r(x)\|_U^p\,\mu(dx)\right]^{\frac{1}{p}}\le c\|R\|_{\gamma(U, L^p(D,\mu))}.
\end{equation*}
\end{proposition}
\begin{proof}
See \cite[Theorem 2.3]{BrezNeerv03}.
\end{proof}

\subsection{Hypotheses and further notation}
\label{sec:hypotheses}
Throughout the rest of the paper, the following is assumed: 
	\begin{itemize}
	\itemsep-0.5em
		\item $(D, \mu)$ is a separable $\sigma$-finite measure space (see Remark \ref{rem:separable}) and the symbol $L^p$ is used to denote the Lebesgue space $L^p(D,\mu)$,
		\item $U$ is a real separable Hilbert space,
		\item $B=(B_t, t\geq 0)$ is an $\alpha$-regular $U$-cylindrical Volterra process (see Definition \ref{def:cVP}).
	\end{itemize}
	
Additionally, the process $B$ is assumed to have the following property:
	\begin{itemize}
	\itemsep-0.5em
	\item There exists $N\in\mathbb{N}_0$, such that $B_t: U\rightarrow \mathcal{H}:=\bigoplus_{i=0}^N\mathcal{H}_i$ for all $t\geq 0$ where $\mathcal{H}_i$ is the $i$-th Wiener chaos (see subsection \ref{sec:chaos}).
	\end{itemize}

\section{Stochastic integration with respect to Volterra processes in \texorpdfstring{$L^p$}{Lp}}
\label{sec:stoch_integration_Lp}

In this section, a stochastic integral $I_T(G)$ with respect to the cylindrical Volterra process $B$ is defined and characterization of integrable operators $G$ is given. 
	
\begin{definition}
Let $T>0$. An operator $G\in\L(U,L^p(D;\mathcal{D}))$ is called elementary if 
	\begin{equation}
	\label{eq:elementary_operator}
		[Gu](x)(t)=\sum_{k=1}^mg_k(t)\langle u,e_k\rangle_Uf_k(x)
	\end{equation}
holds for every $u\in U$, every $t\in [0,T]$, and $\mu$-almost every $x\in D$; where $m\in\mathbb{N}$, $\{e_k\}$ is a complete orthonormal basis of $U$, $\{g_k\}\subset\mathcal{C}^1(0,T)$, and $\{f_k\}\subset L^p$.
\end{definition}

Let $G$ be an elementary operator of the form \eqref{eq:elementary_operator}. Consider the linear operator $I_T$ given by 
	\begin{equation*}
		I_T(G) := \sum_{k=1}^m\left(\int_0^Tg_k(r)\d{b}_r^{(k)}\right)f_k
	\end{equation*}
where $b^{(k)}=B(e_k)$. As usual, we have to extend the operator $I_T$ to a larger space of integrable functions. The next proposition shows that the natural space of integrands is the space $\gamma(U, L^p(D;\mathcal{D}))$.

\begin{proposition}
\label{prop:charakterization_of_integrands}
Let $1\leq p<\infty$. A bounded linear operator $G:U\rightarrow L^p(D;\mathcal{D})$ is stochastically integrable if and only if $G\in\gamma(U,L^p(D;\mathcal{D}))$. In this case, 
	\begin{equation*}
		\|I_T(G)\|_{L^q(\Omega;L^p)} \eqsim \|G\|_{\gamma(U,L^p(D;\mathcal{D}))}
	\end{equation*}
holds for every $1\leq q<\infty$. 
\end{proposition}

\begin{proof}
Let $G\in\L(U,L^p(D;\mathcal{D}))$ be elementary. Let $\{\delta_k\}$ be independent standard centered Gaussian random variables. Using successively Proposition \ref{prop:radonifying_sum}, Proposition \ref{prop:radon_norm}, the definition of the $L^p(D;\mathcal{D})$, Proposition \ref{prop:HC_in_chaos} (real centered Gaussian random variables belong to $\mathcal{H}_1$) and the independence of $\delta_k$ and $\delta_l$ for $k\neq l$, we obtain that $G$ is $\gamma$-radonifying and 
	\begin{align*}
	\label{eq:radonifying_norm}
		\|G\|_{\gamma(U,L^p(D;\mathcal{D}))}^p & \eqsim \E\left\|\sum_{k=1}^m\delta_kg_kf_k\right\|^p_{L^p(D;\mathcal{D})}\\
		& = \int_D\E\left\|\sum_{k=1}^m\delta_kg_kf_k(x)\right\|_\mathcal{D}^p\mu(\d{x})\\
		& \eqsim \int_D\left(\E\left\|\sum_{k=1}^m\delta_kg_kf_k(x)\right\|_\mathcal{D}^2\right)^\frac{p}{2}\mu(\d{x})\\
		& \eqsim \int_D\left(\sum_{k=1}^m\|g_k\|_{\mathcal{D}}^2|f_k(x)|^2\right)^\frac{p}{2}\mu(\d{x}).
	\end{align*}
On the other hand, using the definition of $I_T(G)$ and the definition of the $L^p$ norm, we obtain
	\begin{equation*}
		\E\left\|I_T(G)\right\|^p_{L^p} = \E\left\|\sum_{k=1}^m\left(\int_0^Tg_k(r)\d{b}_r^{(k)}\right)f_k\right\|^p_{L^p}  = \int_D\E\left|\sum_{k=1}^m\left(\int_0^Tg_k(r)\d{b}_r^{(k)}\right)f_k(x)\right|^p\mu(\d{x}).
	\end{equation*}
Let us denote
	\begin{equation*}
		\varepsilon_k := \int_0^Tg_k(r)\d{b}_r^{(k)}, \quad \quad I(x):=\E\left|\sum_{k=1}^m\varepsilon_kf_k(x)\right|^p.
	\end{equation*}
Since $b_t\in\mathcal{H}$ for all $t\geq 0$, we have that every one-dimensional elementary integral of the form $\sum_{k=0}^{n-1}\varphi_i(b_{t_{i+1}}-b_{t_i})\in\mathcal{H}$. If $\varphi\in\mathcal{D}$ and $\{\varphi_k\}$ is a sequence of step functions such that $\varphi_k\rightarrow\varphi$ in $\mathcal{D}$, then $\mathcal{H}\ni i_T(\varphi_k)\rightarrow i_T(\varphi)$ in $L^2(\Omega)$ and hence, $i_T(\varphi)\in\mathcal{H}$. This means, that $\varepsilon_k\in\mathcal{H}$ for every $k=1,\ldots, m$ and also $\sum_{k=1}^m\varepsilon_kf_k(x)\in\mathcal{H}$. Hence, we obtain
	\begin{equation*}
		I(x)  = \E\left|\sum_{k=1}^m\varepsilon_kf_k(x)\right|^p
		\eqsim\left(\E\left|\sum_{k=1}^m\varepsilon_kf_k(x)\right|^2\right)^\frac{p}{2}
	\end{equation*}
by Proposition \ref{prop:HC_in_chaos}. Using the fact that $\varepsilon_k$ and $\varepsilon_l$ for $k\neq l$ are uncorrelated (since $b^{(k)}$ and $b^{(l)}$ are uncorrelated for $k\neq l$), it follows that 
	\begin{equation*}
		\E\|I_T(G)\|^p_{L^p}  \eqsim \int_D\left(\sum_{k=1}^m\E(\varepsilon_k^2)\, |f_k(x)|^2\right)^\frac{p}{2}\mu(\d{x}) = \int_D\left(\sum_{k=1}^m\|g_k\|_\mathcal{D}^2|f_k(x)|^2\right)^\frac{p}{2}\mu(\d{x})
	\end{equation*}
by \eqref{eq:Itoisometry}. Proposition \ref{prop:HC_in_chaos} yields
	\begin{equation}
	\label{eq:norm_of_elementary_integral}
		\|I_T(G)\|_{L^q(\Omega;L^p)} \eqsim \|I_T(G)\|_{L^p(\Omega;L^p)} \eqsim \|G\|_{\gamma(U,L^p(D;\mathcal{D}))}
	\end{equation}
for any $1\leq q<\infty$. Let now $G\in \gamma(U,L^p(D;\mathcal{D}))$ be arbitrary and let $\{G_k\}$ be a sequence of elementary operators such that $G_k\rightarrow G$ in $\gamma(U,L^p(D;\mathcal{D}))$. By \eqref{eq:norm_of_elementary_integral}, we have that
	\begin{equation*}
		\|I_T(G_k)-I_T(G_l)\|_{L^2(\Omega;L^p)} = \|I_T(G_k-G_l)\|_{L^2(\Omega;L^p)}\eqsim\|G_k-G_l\|_{\gamma(U,L^p(D;\mathcal{D}))}
	\end{equation*}
which tends to zero as $k,l\rightarrow\infty$. Hence, $\{I_T(G_k)\}$ is a Cauchy sequence in $L^2(\Omega;L^p)$ and since this is a Banach space, there must be a limit $Z$. The stochastic integral of $G$ is then defined as the map $I_T: G\mapsto Z$. Applying Proposition \ref{prop:HC_in_chaos} again yields the claim.
\end{proof}

\begin{remark}
\label{rem:properties_of_integral}
The stochastic integral of $G\in\gamma(U,L^p(D;\mathcal{D}))$ with respect to $B$ can be written as
	\begin{equation*}
		\int_0^TG(r)\d{B}_r :=  I_T(G) = \sum_{n} \int_0^TGe_n\d{b}_r^{(n)}.
	\end{equation*}
\end{remark}

\begin{corollary}
\label{cor:bound_elem_Lp_int}
Let $\frac{2}{1+2\alpha} \leq p<\infty$. The space $L^\frac{2}{1+2\alpha}([0,T];\gamma(U,L^p))$ is continuously embedded in $\gamma(U,L^p(D;\mathcal{D}))$ and 
	\begin{equation}
	\label{eq:bound_on_integral_of_G}
		\|I_T(G)\|_{L^q(\Omega;L^p)} \lesssim \|G\|_{L^\frac{2}{1+2\alpha}([0,T];\gamma(U,L^p))}
	\end{equation}
holds for every $1\leq q<\infty$ and $G\in L^\frac{2}{1+2\alpha}([0,T];\gamma(U, L^p))$.
\end{corollary}

\begin{proof}
Let $G$ be elementary. From the proof of Proposition \ref{prop:charakterization_of_integrands} we have that 
	\begin{equation*}
		\|I_T(G)\|_{L^q(\Omega;L^p)} \eqsim \left(\int_D\left(\sum_{k=1}^m \|g_k\|_\mathcal{D}^2|f_k(x)|^2\right)^\frac{p}{2}\mu(\d{x})\right)^\frac{1}{p}.
	\end{equation*}
By \eqref{eq:estimate_norm_of_elementary_integral} and the Cauchy-Schwarz inequality it follows that
	\begin{equation}
	\label{eq:est_on_kernel_sum}
		\sum_{k=1}^m\|g_k\|_\mathcal{D}^2|f_k(x)|^2 \lesssim \int_0^T\int_0^T\|J(u,x)\|_U\|J(v,x)\|_U |u-v|^{2\alpha-1}\d{u}\d{v}
	\end{equation}
where 
	\begin{equation*}
		J(u,x):=\sum_{k=1}^mg_k(u)f_k(x)e_k.
	\end{equation*}
Using \eqref{eq:est_on_kernel_sum}, the Hardy-Littlewood inequality and the Minkowski inequality successively yields
	\begin{align*}
		\|I_T(G)\|_{L^q(\Omega;L^p)} & \lesssim \left(\int_D\left(\int_0^T\int_0^T\|J(u,x)\|_U\|J(v,x)\|_U|u-v|^{2\alpha-1}\d{u}\d{v}\right)^\frac{p}{2}\mu(\d{x})\right)^\frac{1}{p}\\
		& \lesssim \left(\int_D\left(\int_0^T\|J(u,x)\|_U^\frac{2}{1+2\alpha}\d{u}\right)^{p\left(\alpha+\frac{1}{2}\right)}\mu(\d{x})\right)^\frac{1}{p}\\
		& \lesssim \left(\int_0^T\left(\int_D\|J(u,x)\|_U^p\mu(\d{x})\right)^\frac{2}{p(1+2\alpha)}\d{u}\right)^{\alpha+\frac{1}{2}}\\
		& \eqsim \left(\int_0^T\|G(u)\|_{\gamma(U,L^p)}^\frac{2}{1+2\alpha}\d{u}\right)^{\alpha+\frac{1}{2}}.
	\end{align*}
The claim for $G\in L^\frac{2}{1+2\alpha}([0,T];\gamma(U,L^p))$ follows by a standard approximation argument.
\end{proof}

\begin{remark}
Let us mention that in the case of fBm with $H>1/2$ or the Rosenblatt process, the condition $\frac{2}{1+2\alpha}\leq p<\infty$ reads as $1\leq pH<\infty$.
\end{remark}

\section{Stochastic evolution equation in \texorpdfstring{$L^p$}{Lp}}
In the rest of the paper, we assume that $\frac{2}{1+2\alpha}\leq p<\infty$. Consider the following stochastic differential equation
	\begin{equation}
	\label{eq:SDE}
		\left\{\begin{array}{ccl}
					\d X_t & = & AX_t\d t + \varPhi\d B_t, \quad t\geq 0,\\
						X_0 & = & x,
				\end{array}\right.
	\end{equation}
where $x\in L^p$, $A:\Dom(A)\rightarrow L^p$, $\Dom(A)\subset L^p$, is an infinitesimal generator of an analytic, strongly continuous semigroup of linear operators $(S(t), t\geq 0)$ acting on $L^p$, $\varPhi\in\L(U,L^p)$ and $B=(B_t, t\geq 0)$ is a $U$-cylindrical Volterra process satisfying the hypotheses from subsection \ref{sec:hypotheses}. The solution to \eqref{eq:SDE} is considered in the mild form, i.e.
	\begin{equation}
	\label{eq:mild_solution}
		X_t = S(t)x+\int_0^tS(t-r)\varPhi\d B_r, \quad t\geq 0.
	\end{equation}

\begin{proposition}
\label{prop:meas_paths}
Assume that $S(u)\varPhi\in \gamma(U,L^p)$ for every $u>0$ and that there is a $T_0>0$ such that 
	\begin{equation}
	\label{eq:assumption_on_existence}
		\int_0^{T_0} \|S(r)\varPhi\|_{\gamma(U,L^p)}^\frac{2}{1+2\alpha}\d{r} <\infty.
	\end{equation}
Then, the solution $X$, given by \eqref{eq:mild_solution}, is well-defined, $L^p$-valued, and mean-square right continuous. In particular, $X$ admits a version with measurable sample paths. 
\end{proposition}

\begin{proof}
\underline{Existence:} By the same arguments as in the proof pro Corollary \ref{cor:bound_elem_Lp_int}, we need to show that the following is finite for every $t>0$:
	\begin{equation*}
		\E\left\|\int_0^tS(t-r)\varPhi\d{B_r}\right\|^2_{L^p}\lesssim \left(\int_0^t\|S(t-r)\varPhi\|_{\gamma(U,L^p)}^\frac{2}{1+2\alpha}\d{r}\right)^{1+2\alpha} = \left(\int_0^t\|S(r)\varPhi\|_{\gamma(U,L^p)}^\frac{2}{1+2\alpha}\d{r}\right)^{1+2\alpha}.
	\end{equation*}
First assume that $t\in[0,T_0]$. Then, we have that 
	\begin{equation*}
		\int_0^t\|S(r)\varPhi\|_{\gamma(U,L^p)}^\frac{2}{1+2\alpha}\d{r} \leq \int_0^{T_0}\|S(r)\varPhi\|_{\gamma(U,L^p)}^\frac{2}{1+2\alpha}\d{r} <\infty.
	\end{equation*}
If $t\in (T_0,\infty)$, then we can write
	\begin{equation*}
		\int_0^t\|S(r)\varPhi\|_{\gamma(U,L^p)}^\frac{2}{1+2\alpha}\d{r} = \int_0^{T_0}\|S(r)\varPhi\|_{\gamma(U,L^p)}^\frac{2}{1+2\alpha}\d{r} + \int_{T_0}^t\|S(r)\varPhi\|_{\gamma(U,L^p)}^\frac{2}{1+2\alpha}\d{r}.
	\end{equation*}
The last integral is finite for every $t\in (T_0,\infty)$ by the following semigroup property: If $\omega>0$ is such that $\|S(u)\|_{\mathcal{L}(L^p)}\lesssim e^{\omega u}$ for $u>0$, then for $r>T_0$ we have that
	\begin{equation}
	\label{eq:semigroup_for_large_t}
		\|S(r)\varPhi\|_{\gamma(U,L^p)} = \|S(r-T_0)S(T_0)\varPhi\|_{\gamma(U,L^p)} \lesssim e^{\omega(r-T_0)}\|S(T_0)\varPhi\|_{\gamma(U,L^p)}.
	\end{equation}

\underline{Mean-square right continuity:} Notice first that we can write
	\begin{align*}
		\left\|\int_0^tS(t-r)\varPhi\d{B}_r - \int_0^sS(s-r)\varPhi\d{B}_r\right\|_{L^2(\Omega;L^p)} & = \\
		& \hspace{-6cm} = \left\|\int_s^tS(t-r)\varPhi\d{B}_r + \int_0^s(S(t-s)-I)S(s-r)\varPhi\d{B}_r\right\|_{L^2(\Omega;L^p)}\\
		& \hspace{-6cm} \leq \left(\E\left\|\int_s^tS(t-r)\varPhi\d{B}_r\right\|^2_{L^p}\right)^\frac{1}{2} + \left(\E\left\|\int_0^s(S(t-s)-I)S(s-r)\varPhi\d{B}_r\right\|^2_{L^p}\right)^\frac{1}{2}.
	\end{align*}
for $0<s<t$. For the first integral, write
	\begin{align*}
		I_1(s,t):=\E\left\|\int_s^tS(t-r)\varPhi\d{B}_r\right\|^2_{L^p} & \lesssim \left(\int_s^t\|S(t-r)\varPhi\|_{\gamma(U,L^p)}^\frac{2}{1+2\alpha}\d{r}\right)^{1+2\alpha}\\
		& = \left(\int_0^{t-s}\|S(r)\varPhi\|_{\gamma(U,L^p)}^\frac{2}{1+2\alpha}\d{r}\right)^{1+2\alpha}
	\end{align*}
Hence, $I_1(s,t)$ tends to $0$ as $t\rightarrow s+$. The second integral can be estimated as follows:
	\begin{align*}
		I_2(s,t)& :=\E\left\|\int_0^s(S(t-s)-I)S(s-r)\varPhi\d{B}_r\right\|^2_{L^p} \lesssim\\
		& \hspace{3cm} \lesssim \left(\int_0^s\|(S(t-s)-I)S(s-r)\varPhi\|_{\gamma(U,L^p)}^\frac{2}{1+2\alpha}\d{r}\right)^{1+2\alpha}.
	\end{align*}
By the definition of the $\gamma$-radonifying norm and Proposition \ref{prop:radonifying_sum}, if we take a centered sequence $\{\delta_n\}$ of independent, standard Gaussian random variables, we have that
	\begin{equation*}
		\|(S(t-s)-I)S(s-r)\varPhi\|_{\gamma(U,L^p)}  = \left(\E\left\|\sum_{n=1}^\infty\delta_n(S(t-s)-I)S(s-r)\varPhi e_n\right\|^2_{L^p}\right)^\frac{1}{2}
	\end{equation*}
Strong continuity of the semigroup $(S(t), t\geq 0)$ implies that, for every $r\in [0,s)$,
	\begin{equation*}
		\left\|(S(t-s)-I)\left(\sum_{n=1}^\infty \delta_n S(s-r)\varPhi e_n\right)\right\|_{L^p}^2\rightarrow 0 \quad \Pas
	\end{equation*}
as $t\rightarrow s+$. Moreover, we have that
	\begin{equation*}
		\left\|(S(t-s)-I)\left(\sum_{n=1}^\infty\delta_n S(s-r)\varPhi e_n\right)\right\|_{L^p}^2\lesssim\left\|\sum_{n=1}^\infty\delta_nS(s-r)\varPhi e_n\right\|_{L^p}^2 \quad \Pas
	\end{equation*} 
and
	\begin{equation*}
		\E\left\| \sum_{n=1}^\infty\delta_nS(s-r)\varPhi e_n\right\|_{L^p}^2 = \|S(s-r)\varPhi\|_{\gamma(U,L^p)}^2 <\infty
	\end{equation*}
by the assumptions of the proposition. Therefore, we obtain
	\begin{equation*}
		\| (S(t-s)-I)S(s-r)\varPhi\|_{\gamma(U,L^p)}^\frac{2}{1+2\alpha}\rightarrow 0
	\end{equation*}
as $t\rightarrow s+$ for every $r\in [0,s)$ by the Lebesgue Dominated Convergence Theorem (DCT). Furthermore, we have that
	\begin{equation*}
		\|(S(t-s)-I)S(s-r)\varPhi\|_{\gamma(U,L^p)}\leq \|S(t-s)-I\|_{\L(L^p)}\|S(s-r)\varPhi\|_{\gamma(U,L^p)}
	\end{equation*}
and 
	\begin{equation*}
		\int_0^s\|S(s-r)\varPhi\|_{\gamma(U, L^p)}^\frac{2}{1+2\alpha}\d{r} <\infty
	\end{equation*}
by the first part of the proof. This yields $I_2(s,t)\rightarrow 0$ as $t\rightarrow s+$ by the DCT. Hence, we have proved mean-square right continuity of the process $X$. The existence of a version with measurable sample paths follows by standard arguments (see e.g. \cite[Proposition 3.2]{DaPratoZab14book}).
\end{proof}

Continuity of the solution $X$ to \eqref{eq:SDE} is discussed now. Since the semigroup $S$ is analytic, there is $\lambda\in\R$ such that the operator $(\lambda I- A)$ is strictly positive. Let us thus denote, for $\delta\geq 0$, 
	\begin{equation*}
		V_{\delta,p} :=\Dom((\lambda I-A)^\delta)\subset L^p.
	\end{equation*}
Equipped with the graph norm topology, the space $V_{\delta,p}$ is a Banach space. The main result follows.

\begin{proposition}
\label{prop:V_dp-continuity}
Assume that $S(u)\varPhi\in\gamma(U,L^p)$ for all $u>0$ and that there are $T_0>0$, $\delta\geq 0$ and $\beta>0$ such that $x\in V_{\delta,p}$,
\begin{equation*}
\beta +\delta <\alpha+\frac{1}{2}
\end{equation*}
 and
\begin{equation}
	\label{eq:Vd_cts}
		\int_0^{T_0}\left(r^{-(\beta+\delta)}\|S(r)\varPhi\|_{\gamma(U,L^p)}\right)^\frac{2}{1+2\alpha}\d{r} <\infty.
	\end{equation}
Then, $X$ has a version which belongs to $\mathcal{C}^{\nu}([0,T];V_{\delta,p})$ a.s. for every $\nu\in [0,\beta)$.
\end{proposition}

\begin{proof}
\underline{Step 1:} We show that the integral
	\begin{equation*}
		\int_0^t(\lambda I-A)^\delta S(t-r)\varPhi\d B_r
	\end{equation*}
exists for all $t>0$. First we have to notice that since $S(u)\varPhi\in\gamma(U,L^p)$ for each $u>0$, it follows that also $(\lambda I-A)^\delta S(u)\varPhi\in\gamma(U,L^p)$ for all $u>0$. In the same way as in Corollary \ref{cor:bound_elem_Lp_int}, we can obtain
	\begin{align*}
		\E\left\|\int_0^t(\lambda I-A)^\delta S(t-r)\varPhi\d B_r\right\|_{L^p}^2 & \lesssim \left(\int_0^t\|(\lambda I-A)^\delta S(t-r)\varPhi\|_{\gamma(U,L^p)}^\frac{2}{1+2\alpha}\d r\right)^{1+2\alpha}\\
		& \lesssim \left(\int_0^{\frac{t}{2}}\left(r^{-\delta}\|S(r)\varPhi\|_{\gamma(U,L^p)}\right)^\frac{2}{1+2\alpha}\d r\right)^{1+2\alpha}
	\end{align*}
since the integrand can be estimated by
	\begin{align*}
		\|(\lambda I-A)^\delta S(u)\varPhi\|_{\gamma(U,L^p)} &\leq \left\|(\lambda I-A)^\delta S\left(\frac{u}{2}\right)\right\|_{\L( L^p)}\left\|S\left(\frac{u}{2}\right)\varPhi\right\|_{\gamma(U,L^p)} \\
		& \lesssim u^{-\delta}\left\|S\left(\frac{u}{2}\right)\varPhi\right\|_{\gamma(U,L^p)}
	\end{align*}
for $u>0$. Now, if $t<2T_0$, we can only enlarge the integration bounds and use \eqref{eq:Vd_cts}. For $t>2T_0$, we use the semigroup property \eqref{eq:semigroup_for_large_t}.
	
\underline{Step 2:} Since $(\lambda I-A)^\delta$ is a closed operator, it follows that $\int_0^tS(t-r)\varPhi\d B_r\in V_{\delta,p}$. Moreover,
	\begin{equation*}
		(\lambda I-A)^\delta\int_0^tS(t-r)\varPhi\d B_r = \int_0^t(\lambda I-A)^\delta S(t-r)\varPhi\d B_r \quad \Pas
	\end{equation*}
for every $t>0$. 

\underline{Step 3:} We use the factorization technique to show that the solution admits a $V_{\delta,p}$-valued version with H\"{o}lder continuous sample paths if $x\in V_{\delta,p}$. Fix $T>0$. Similarly as above, we can show that
	\begin{equation*}
		\E\left\|\int_0^u(u-r)^{-\beta}(\lambda I-A)^\delta S(u-r)\varPhi\d B_r\right\|_{L^p}^2\lesssim \left(\int_0^{u/2} \left(r^{-(\beta+\delta)}\|S(r)\varPhi\|_{\gamma(U, L^p)}\right)^\frac{2}{1+2\alpha}\d r\right)^{1+2\alpha}
	\end{equation*}
and the assumption \eqref{eq:Vd_cts} assures that this is finite for all $u>0$. Since again, $(\lambda I-A)^\delta$ is closed, we have that 
	\begin{equation*}
		\int_0^u(u-r)^{-\beta}S(u-r)\varPhi\d B_r\,\,\in\,\, V_{\delta,p}\quad\Pas
	\end{equation*}
for all $u>0$ and the integral commutes with $(\lambda I-A)^\delta$. This means that we can define an $L^p$-valued process $Y^\delta=(Y^\delta_u, u\in [0,T])$ by 
	\begin{equation}
	\label{eq:Y_in_factorization}
		Y^\delta_u:=\int_0^u(u-r)^{-\beta}(\lambda I-A)^\delta S(u-r)\varPhi\d{B}_r.
	\end{equation}
$Y^\delta$ has a version with measurable sample paths which can be shown similarly as in the proof of Proposition \ref{prop:meas_paths} using \eqref{eq:Vd_cts}. Moreover, since 
	\begin{equation*}
		\sup_{t\in (0,T]}\|Y^\delta_t\|_{L^q(\Omega;L^p)} \lesssim\sup_{t\in (0,T]}\left(\int_0^t\|r^{-(\beta+\delta)}S(r)\varPhi\|_{\gamma(U,L^p)}^\frac{2}{1+2\alpha}\d{r}\right)^{\alpha+\frac{1}{2}}<\infty
	\end{equation*}
for every $1\leq q<\infty$ by \eqref{eq:bound_on_integral_of_G} and \eqref{eq:Vd_cts}, we infer that $Y^\delta\in L^q([0,T];L^p)$ $\mathbb{P}$-almost surely for every $1\leq q<\infty$.

Recall that 
	\begin{equation*}
		\int_r^t(t-u)^{\beta-1}(u-r)^{-\beta}\d u=\frac{\pi}{\sin \pi\beta} =:\frac{1}{\Lambda}
	\end{equation*}
holds for $\beta\in (0,1)$ and $r\in [0,t]$. Using this fact, we can write
	\begin{align*}
		\int_0^t (\lambda I-A)^\delta S(t-r)\varPhi\d B_r & = \\
		& \hspace{-3cm} = \Lambda\sum_{n=1}^\infty\int_0^t\int_r^t(t-u)^{\beta-1}(u-r)^{-\beta}(\lambda I-A)^\delta S(t-u)S(u-r)\varPhi e_n\d u\d b^{(n)}_r\\
		& \hspace{-3cm} = \Lambda\sum_{n=1}^\infty \int_0^t(t-u)^{\beta-1}S(t-u)\underbrace{\left(\int_0^u(u-r)^{-\beta}(\lambda I-A)^{\delta}S(u-r)\varPhi e_n\d b_r^{(n)}\right)}_{=:Y_n(u)}\d u\\
		& \hspace{-3cm} = \Lambda\lim_{N\rightarrow\infty} \int_0^t(t-u)^{\beta-1}S(t-u)\underbrace{\left(\sum_{n=1}^NY_n(u)\right)}_{=:Y^N(u)}\d u.
	\end{align*}
The interchange of the order of integration is possible due to the fact that the function 
	\begin{equation*}
		f(u,r) := \bm{1}_{(r,t]}(u)(t-u)^{\beta-1}(u-r)^{-\beta}(\lambda I-A)^\delta S(t-r)\varPhi z
	\end{equation*}
(here $z\in U$) belongs to the mixed Lebesgue space $L^{1,\frac{2}{1+2\alpha}}([0,t];L^p)$ and $f^*(u,r):= f(r,u)$ belongs to $L^{\frac{2}{1+2\alpha},1}([0,t];L^p)$ since $0<\beta+\delta<\alpha+\frac{1}{2}$ (see \cite{BenedekPanzone}). Therefore, a suitable stochastic Fubini theorem can be proved for this particular function (see \cite[Lemma 4.4]{CouMas16} for a similar result). 
In order to interchange the limit $N\rightarrow\infty$ and the integral, consider the following:
	\begin{multline*}
		\left\|\int_0^t(t-u)^{\beta-1}S(t-u)Y^N(u)\d{u} - \int_0^t(t-u)^{\beta-1}S(t-u)Y_u^\delta\d{u}\right\|_{L^1(\Omega;L^p)} \lesssim \\ \lesssim \int_0^t(t-u)^{\beta-1}\|S(t-u)\|_{\L(L^p)}\|Y^N(u)-Y^\delta_u\|_{L^p(\Omega;L^p)}\d{u}.
	\end{multline*}
The norm inside the last integral can be estimated by
	\begin{align*}
		\sup_{N}\sup_{u\in (0,t]} \|Y^N(u)-Y_u^\delta\|_{L^p(\Omega;L^p)}  & \lesssim \\
			& \hspace{-4cm} \lesssim \sup_{N}\sup_{u\in (0,t]} \left(\int_0^ur^{\frac{-2(\beta+\delta)}{1+2\alpha}}\left(\int_D\left(\sum_{n=N+1}^\infty\left|[S(r)\varPhi e_n](x)\right|^2\right)^\frac{p}{2}\mu(\d{x})\right)^\frac{2}{p(1+2\alpha)}\d{r}\right)^{\alpha+\frac{1}{2}} \\
			&  \hspace{-4cm} \lesssim \left(\int_0^t\left(r^{-(\beta+\delta)}\|S(r)\varPhi\|_{\gamma(U,L^p)}\right)^\frac{2}{1+2\alpha}\d{r}\right)^{\alpha+\frac{1}{2}}
	\end{align*}
by similar computations as in the proof of Corollary \ref{cor:bound_elem_Lp_int}. Hence, we can infer by DCT that 
	\begin{equation*}
		\int_0^t (\lambda I-A)^\delta S(t-r)\varPhi\d B_r = \Lambda \int_0^t(t-u)^{\beta-1}S(t-u)Y_u^\delta\d u, \quad \Pas
	\end{equation*}
holds for every $t>0$. Consider the operator 
	\begin{equation*}
		R_{\beta,T}(Z)(t) := \int_0^t(t-u)^{\beta-1}S(t-u)Z(u)\d u, \quad t\in [0,T].
	\end{equation*}
By \cite[Theorem 5.14, (ii)]{DaPratoZab14book}, the operator $R_{\beta, T}$ is bounded from $L^q([0,T];L^p)$ to $\mathcal{C}^{\nu}([0,T];L^p)$ for every $\nu\in [0,\beta-\frac{1}{q})$. Taking $q$ sufficiently large yields the claim of the proposition.
\end{proof}

\begin{remark}
The proof of Proposition \ref{prop:V_dp-continuity} is based on the factorization method (see e.g. \cite{DaPKwaZab87}). In the literature on stochastic convolution in Banach spaces, $\gamma$-boundedness combined with estimates on analytic semigroups is also used (see e.g. \cite{NeeVerWei15,VerWei11}).
\end{remark}

\begin{corollary}
\label{cor: gamma}
Assume that $S(u)\varPhi\in\gamma(U,L^p)$ and that there is a $\gamma\in[0,\alpha+\frac{1}{2})$ such that 
	\begin{equation}
	\label{eq:gamma_estimate}
		\|S(u)\varPhi\|_{\gamma(U, L^p)}\lesssim u^{-\gamma}
	\end{equation} 
for all $u>0$. Then, $X$ has a version which belongs to $\mathcal{C}^\nu([0,T];V_{\delta,p})$ a.s. for every $\nu,\delta\geq 0$ such that $x\in V_{\delta,p}$ and
	\begin{equation*}
		\nu+\delta < \alpha+\frac{1}{2}-\gamma.
	\end{equation*}
\end{corollary}

\begin{proof}
Let $\delta\in [0,\alpha+\frac{1}{2}-\gamma)$ be arbitrary but fixed. Now we can choose $\beta>0$ such that
	\begin{equation*}
		0< \beta +\delta +\gamma < \alpha +\frac{1}{2}.
	\end{equation*}
Then \eqref{eq:Vd_cts} holds since 
	\begin{equation*}
		\int_0^T\left(r^{-(\beta+\delta)}\|S(r)\varPhi\|_{\gamma(U,L^p)}\right)^\frac{2}{1+2\alpha}\d r \lesssim \int_0^Tr^{\frac{-2}{1+2\alpha}(\beta+\delta+\gamma)}\d r <\infty.
	\end{equation*}
Hence, $X\in \mathcal{C}^{\nu}([0,T];V_{\delta,p})$ for every $\nu\in [0,\beta)$ by Proposition \ref{prop:V_dp-continuity}. Taking the supremum over all such $\beta$'s yields the claim.
\end{proof}

\begin{remark}
For the fBm or the Rosenblatt process ($H>\frac{1}{2}$ for both), we obtain that if there is $\gamma\in [0,H)$ such that $\|S(u)\varPhi\|_{\gamma(U,L^p)}\lesssim u^{-\gamma}$ for all $u>0$, then for every $0\leq \delta+\nu<H-\gamma$ with $x\in V_{\delta,p}$, the solution has a version in $\mathcal{C}^\nu([0,T];V_{\delta,p})$.
\end{remark}

\section{Examples}

\subsection{Parabolic equations with pointwise Volterra noise} Consider the following parabolic equation
	\begin{equation*}
		\partial_tu = \Delta u + \delta_z\eta, \quad\mbox{on}\quad {\R}_+\times\O
	\end{equation*}
with the initial condition $u(0,\cdot)=f$ on $\O$ and the Dirichlet boundary condition $u|_{{\R}_+\times\partial\O}=0$. Given a point $z\in\O$, $\O\subset\R^d$ an open bounded domain with $\mathcal{C}^1$ boundary $\partial\O$, $\delta_z$ denotes the Dirac distribution at $z$. 

This formal system can be rewritten as the stochastic evolution equation \eqref{eq:SDE}. The noise process $\eta$ is the formal derivative of a scalar $\alpha$-regular Volterra process $b=(b_t, t\geq 0)$ which belongs to a finite Wiener chaos $\mathcal{H}$. We assume that $p\geq \frac{2}{1+2\alpha}$ and take $x:=f\in L^p(\O)$; $A:=\Delta|_{\Dom(A)}$ with 
	\begin{equation*}
		\Dom(A):=W^{2,p}(\O)\cap W_0^{1,p}(\O),
	\end{equation*}
which is a generator of an analytic semigroup on $L^p(\mathcal{O})$; and $\varPhi$ which is given by $\varPhi a = a\delta_z$ for $a\in\R=:U$. By the Sobolev embedding (see e.g. \cite[Theorem 8.2]{NezPalVal12}), for every $\varepsilon\in (0, 1-\frac{d}{2p})$ we have that $\varPhi\in\L(\R,V_{\varepsilon-1,p})$ since $V_{\delta,p}\subset W^{2\delta,p}(\O)$. 

Note that
	\begin{equation*}
		\|S(r)\varPhi\|_{\gamma(\R,L^p)} \leq \|S(r)\|_{\L(V_{\varepsilon-1,p}, L^p)}\|\varPhi\|_{\L(\R,V_{\varepsilon-1},p)} \lesssim r^{\varepsilon-1}
	\end{equation*}
for $r>0$. Thus, we can apply Corollary \ref{cor: gamma} with $\gamma := 1-\varepsilon$. 

If 
	\begin{equation}
	\label{eq:alpha_pointwise}
		\alpha>\frac{d}{2p}-\frac{1}{2},
	\end{equation}
then we can choose $\varepsilon$ such that
	\begin{equation*}
		\varepsilon\in\left(\frac{1}{2}-\alpha, 1-\frac{d}{2p}\right)
	\end{equation*}	
so that, by Corollary \ref{cor: gamma}, there is a version of the solution $X$ which belongs to $\mathcal{C}^\nu([0,T];V_{\delta,p})$ for every $\nu,\delta\geq 0$ such that $x\in V_{\delta,p}$ and $\nu+\delta <\alpha+\varepsilon-\frac{1}{2}$. Taking the supremum over all such $\varepsilon$ yields that $X$ has a version $\tilde{X}$ such that 
	\begin{equation*}
\tilde{X}\in\mathcal{C}^\nu([0,T];V_{\delta,p}), \quad\mbox{ for all }\quad  \nu+\delta <\alpha+\frac{1}{2} -\frac{d}{2p}.
	\end{equation*}
Note that \eqref{eq:alpha_pointwise} does not pose additional constraints on $\alpha$ if $d=1$; it excludes the case $p=\frac{2}{1+2\alpha}$ if $d=2$ and; finally, for $d\geq 3$, \eqref{eq:alpha_pointwise} can only be satisfied if $2p>d$.

If, however, the stronger condition 
	\begin{equation}
	\label{eq:alpha_pointwise_2}
		\alpha > \frac{d}{p}-\frac{1}{2}
	\end{equation}
is satisfied, then, by the Sobolev embedding, we have that $X$ has a version $\tilde{X}$ such that
	\begin{equation}
	\label{eq:X_in_C}
		\tilde{X} \in \mathcal{C}^\nu([0,T];\mathcal{C}^{2\delta-\frac{d}{p}}(\bar{\O}))
	\end{equation}	
for every $\nu\geq 0$ and $\delta>\frac{d}{2p}$ such that $x\in V_{\delta,p}$, $\nu+\delta<\alpha +\frac{1}{2}-\frac{d}{2p}$. Note that if $d=1$, the condition \eqref{eq:alpha_pointwise_2} excludes the case $p=\frac{2}{1+2\alpha}$, and, in higher dimensions ($d\geq 2$), it can only be satisfied if $p>d$.

\subsection{Parabolic equations with distributed Volterra noise} Let $m\in\mathbb{N}$ and consider the following parabolic equation
	\begin{equation*}
		\partial_t u = L_{2m}u + \eta\quad \mbox{ on } \quad {\R}_+\times\O
	\end{equation*}
with the initial condition $u(0,\cdot)=f$ which belongs to the space $L^p(\O)$; and the Dirichlet boundary condition
	\begin{equation*}
		\left.\frac{\partial^k u}{\partial\bm{v}^k}\right|_{{\R}_+\times\partial\O}=0
	\end{equation*}
for $k\in \{0,\ldots, m-1\}$ where $\frac{\partial}{\partial\bm{v}}$ denotes the conormal derivative. Here, $\O\subset\R^d$ is an open bounded domain with smooth boundary and $L_{2m}$ is a differential operator of order $2m$, i.e.
	\begin{equation*}
		L_{2m} = \sum_{|k|\leq 2m}a_k(\cdot)\partial^k,
	\end{equation*}
with $a_k\in\mathcal{C}_b^\infty(\O)$ which is assumed to be uniformly elliptic. The considered noise $\eta$ is Volterra in time and can be both white or correlated in space. This system can be rewritten as the stochastic evolution equation \eqref{eq:SDE} in $L^p(\O)$. Indeed, let $U:=L^2(\O)$ and $B=(B_t, t\geq 0)$ be a $U$-cylindrical $\alpha$-regular Volterra process which satisfies the hypotheses of section \ref{sec:hypotheses}. Then the noise $\eta$ is formally given by
	\begin{equation*}
		\eta(t, \cdot) =\varPhi\frac{\d}{\d{t}}B_t
	\end{equation*}
where $\varPhi\in\L(U,L^p(\O))$ determines the space correlation of the noise process $\eta$. Assume that $\frac{2}{1+2\alpha}\leq p <\infty$ and take $x=f\in L^p(\O)$ and $A:=L_{2m}|_{\Dom(A)}$ where
	\begin{equation*}
		\Dom(A) := \left\{f\in W^{2m,p}(\O): \frac{\partial^kf}{\partial\bm{v}^k} = 0 \mbox{ on } \partial\O\mbox{ for } k\in\{0,
		\ldots, m-1\}\right\}.
	\end{equation*}
The operator $A$ generates an analytic semigroup $(S(t), t\geq 0)$ on $L^p(\O)$. By standard estimates on the Green function, we have that
	\begin{equation*}
		\|S(r)\varPhi\|_{\gamma(U,L^p(\O))}\lesssim r^{-\frac{d}{4m}}
	\end{equation*}
for $r>0$ then we can use Corollary \ref{cor: gamma} with $\gamma := \frac{d}{4m}$.

Thus, if 
	\begin{equation}
	\label{eq:alpha_2m_1}
		\alpha > \frac{d}{4m}-\frac{1}{2},
	\end{equation}
then, by Corollary \ref{cor: gamma}, we have that for every $\nu, \delta\geq 0$ such that $x\in V_{\delta, p}$ and $\nu+\delta <\alpha+\frac{1}{2}-\frac{d}{4m}$ the solution $X$ has a version in $\mathcal{C}^\nu([0,T];V_{\delta,p})$. If, moreover, 
	\begin{equation}
	\label{eq:alpha_2m_2}
		\alpha > \frac{d}{4m} -\frac{1}{2} +\frac{d}{2mp},
	\end{equation}
then by the Sobolev embedding, we have that for every $\nu\geq 0$ and $\delta>\frac{d}{2mp}$ such that $x\in V_{\delta,p}$ and $\nu+\delta <\alpha+\frac{1}{2}-\frac{d}{4m}$, the solution has a version in the space $\mathcal{C}^{\nu}([0,T];\mathcal{C}^{2\delta m-\frac{d}{p}}(\bar{\O}))$.

Note that the condition \eqref{eq:alpha_2m_1} can be only satisfied if $d<4m$ and the condition \eqref{eq:alpha_2m_2} can only be satisfied if $d<\frac{4mp}{p+2}$. In the particular case of the stochastic heat equation (i.e. $m=1$), we have that if $d<4$, it is possible to take sufficiently smooth Volterra noise so that the solution is time (H\"{o}lder) continuous in the space $V_{\delta, p}$ and if, moreover, we have that $d < \frac{4p}{p+2}$, then we may even obtain (H\"{o}lder) continuity in the spatial variable. If the initial condition is regular (i.e. $p>6$), the space-time continuity may be obtained in dimensions $d=1,2,3$. 

\bigskip

\textbf{Acknowledgement:} The authors wish to thank the anonymous referee for their careful reading of the paper and for providing useful comments and suggestions.

\end{document}